\documentclass[12pt,reqno]{amsart}

\usepackage{setspace}

\usepackage{caption}

\usepackage
{hyperref, epigraph}

\usepackage{breakurl}   

\usepackage{graphicx} 

\usepackage[export]{adjustbox}  

\usepackage{enumitem} 
\usepackage{breakurl}   

\numberwithin{equation}{section}

\newcommand\R {{\mathbb R}}

\title[When does a hyperbola meet its asymptote?]
{When does a hyperbola meet its asymptote?  Bounded infinities,
  fictions, and contradictions in Leibniz}

\author[M. Katz]{Mikhail G. Katz}\address{M.~Katz, Department of
  Mathematics, Bar Ilan University, Ramat Gan 5290002
  Israel}\email{katzmik@math.biu.ac.il}

\author[D. Sherry]{David Sherry}\address{D. Sherry, Department of
  Philosophy, Northern Arizona University, Flagstaff, AZ 86011,
  US}\email{David.Sherry@nau.edu}

\author[M. Ugaglia]{Monica Ugaglia}\address{M. Ugaglia, Il Gallo
  Silvestre, Localit\`a Collina 38, Montecassiano,
  Italy}\email{monica.ugaglia@gmail.com}

\begin{document}

\thispagestyle{empty}


\keywords{Infinitesimal calculus; useful fiction; infinity;
  infinitesimals; ideal perspective point; Desargues}

\subjclass[2020]{01A45}

\begin{abstract}
In his 1676 text \emph{De Quadratura Arithmetica}, Leibniz
distinguished \emph{infinita terminata} from \emph{infinita
  interminata}.  The text also deals with the notion, originating with
Desargues, of the perspective point of intersection at infinite
distance for parallel lines.  We examine contrasting interpretations
of these notions in the context of Leibniz's analysis of asymptotes
for logarithmic curves and hyperbolas.  We point out difficulties that
arise due to conflating these notions of infinity.  As noted by
Rodr\'{\i}guez Hurtado et al., a significant difference exists between
the Cartesian model of magnitudes and Leibniz's search for a
qualitative model for studying perspective, including ideal points at
infinity.  We show how respecting the distinction between these
notions enables a consistent interpretation thereof.

\emph{Revista Latinoamericana de Filosof\'{\i}a} 49 (2023), no. 2,
241--258.
\end{abstract}

\maketitle

\tableofcontents


\section{Bounded and unbounded infinity}
\label{s1}

A key distinction in Leibniz's approach to geometry and the calculus
is that between bounded infinities (\emph{infinita terminata}) and
unbounded infinities (\emph{infinita interminata}).  As noted by
Knobloch, the distinction was elaborated in Proposition\;11 of the
treatise \emph{De Quadratura Arithmetica} (DQA):
\begin{enumerate}\item[]
[Leibniz] distinguished between two infinites, the bounded infinite
straight line, the recta infinita terminata, and the unbounded
infinite straight line, the recta infinita interminata.  He
investigated this distinction in several studies from the year 1676.
Only the first kind of straight lines can be used in mathematics, as
he underlined in his proof of theorem 11 [i.e., Propositio~XI].%
%
%
\footnote{Knobloch \cite[p.\;97]{Kn99}.}
\end{enumerate}

Leibniz mentioned the distinction in a 29 july\,1698 letter to
Bernoulli.  Here Leibniz analyzes a geometric problem involving
unbounded infinite areas and states an apparent paradox: ``Therefore
the two infinite spaces are equal, and the part is equal to the whole,
which is impossible.''%
\footnote{``Ergo haec duo spatia infinita sibi sunt aequalia, pars
  toti, quod est absurdum'' \cite[p.\;523]{Le98}.}
To resolve the paradox (i.e., the clash with the part-whole
principle), Leibniz exploits the notion of bounded infinity:
\begin{enumerate}\item[]
Properly speaking, the last (ultima) abscissa~$A\,_0B$ is not null, as
if~$O$ were falling on $A$, and the last (ultima)
ordinate~${}_0B\,_0C$ is not unbounded (\emph{interminata}), as
if~${}_0B\,_0C$ were falling on the asymptote.  Rather,~$A\,_0B$ is
infinitely small, and~${}_0B\,_0C$ is infinitely large, but bounded
(\emph{terminata}).%
\footnote{``neque adeo ultimam abscissam~$A\,_0B$ accurate loquendo
  esse nullam, quasi~$O$ incideret in~$A$, nec ultimam ordinatam
  ~${}_0B\,_0C$ esse interminatam, quasi~${}_0B\,_0C$ incideret in
  Asymptotam; sed~$A\,_0B$ esse infinite parvam, et~${}_0B\,_0C$ esse
  infinite magnam, sed terminatam'' (ibid.).  Note that the subscripts
  are on the left in~${}_0B$ and~${}_0C$, as elsewhere in the sequel.}
\end{enumerate}

As noted by Knobloch,%
\footnote{Knobloch \cite[pp.\;267--268]{Kn94}.}
Leibniz also mentioned the distinction in his february 1702
correspondence with Varignon.%
\footnote{Leibniz \cite[p.\;91]{Le02}.}
The distinction enabled him to avoid contradicting the part-whole
principle while still employing infinite magnitudes in analysis and
geometry.

\subsection{Inconsistency of Maxima and Minima}

Leibniz used the term \emph{Maxima} to refer to infinite wholes, and
the term \emph{Minima} to refer to points viewed as constituent parts
of the continuum.  Already in his 1672/3 text ``On Minimum and
Maximum,'' Leibniz rejected both Minima and Maxima in the following
terms:
\begin{enumerate}\item[]
\emph{Scholium}.  We therefore hold that two things are excluded from
the realm of intelligibles: minimum and maximum; the indivisible, or
what is entirely \emph{one}, and \emph{everything}; what lacks parts,
and what cannot be part of another.%
\footnote{Leibniz as translated by Arthur in \cite[p.\,13]{Ar01}.}
\end{enumerate}
Leibniz's rejection of Maxima amounts to the rejection of infinite
wholes (e.g., unbounded lines) as inconsistent, while the rejection of
their counterparts, Minima, amounts to the rejection of putative
simplest constituents of the continuum, i.e., the rejection of a
punctiform continuum.%
\footnote{See \cite{23a}, \cite{23d}.}
To Leibniz, points play only the role of endpoints of line segments.
Thus we find in the Scholium to Proposition~11:



\begin{enumerate}\item[]
The magnitude of an unbounded line, just as that of a point, is beyond
the realm of geometric considerations.%
\footnote{``Nam lineae interminatae magnitudo nullo modo Geometricis
  considerationibus subdita est, non magis quam puncti.''
  \cite[p.\;549]{Le2012}.}
\end{enumerate}

An example of an unbounded infinity is an unending line, or the
continuum.  Such infinities, when taken as a whole, were considered by
Leibniz to lead to a contradiction with the part-whole principle, and
therefore of little use in geometry and calculus.

Unlike unbounded infinity, \emph{bounded infinity} has a pair of
endpoints.  With reference to a line, a Leibnizian bounded infinity
can be thought of as a segment with infinitely separated endpoints
(see Section~\ref{s3}, item\;3).  In the Scholium to Proposition~11,
Leibniz speaks of ``linea terminata qvidem, infinita tamen,''%
\footnote{Leibniz \cite[p.\;549]{Le2012}.}
i.e., a bounded infinite line.

\subsection{Proposition\;11}

In Proposition 11 of DQA, Leibniz discusses bounded infinities and
explains their reciprocal relation to infinitesimals.  He uses the
notation~$(\mu)\mu$ for an infinitesimal.  Here~$\mu$ can be thought
of as the origin, and $(\mu)$ as an infinitely close point.  Leibniz
notes that in order to show that a figure of infinite length bounded
by a curve~$\lambda$ may have a finite magnitude one must proceed as
follows:
\begin{enumerate}\item[]
%
Substitute for a line $\mu\lambda$ [the asymptote] a line
$(\mu)\lambda$, the point $(\mu)$ being taken just above $\mu$ and the
interval $(\mu)\mu$ being infinitely small, so that the ordinate
$(\mu)\lambda$ will be of infinite length.%
\footnote{``pro recta $\mu\lambda$ ponatur recta $(\mu)\lambda$ puncto
  $(\mu)$ paulo supra punctum $\mu$ sumto, intervallo $(\mu)\mu$
  infinite parvo, ita ordinata $(\mu)\lambda$ erit longitudine
  infinita'' \cite[p.\;547]{Le2012}.}
\end{enumerate}
Here Leibniz uses the notation $x\lambda$ for the ordinate of the
point on the curve $\lambda$ with abscissa $x$.  The bounded infinity,
denoted~$(\mu)\lambda$, is the ordinate of the point on the
curve~$\lambda$ corresponding to an infinitesimal abscissa~$(\mu)$.
Leibniz goes on to emphasize that
%
%
``$(\mu)\lambda$ will not be the asymptote.''%
\footnote{%
%
%
``Proinde $(\mu)\lambda$ non erit curvae $D\delta$ asymptotos''
  \cite[p.\;547]{Le2012}.}
Thus bounded infinity is distinct from the asymptote.  Note that the
line~$(\mu)\lambda$ is a subline of the infinite unbounded line with
abscissa~$(\mu)$.

\begin{figure}
\begin{center}
\includegraphics
[scale=0.5,
trim = 0 0 0 0,clip]
{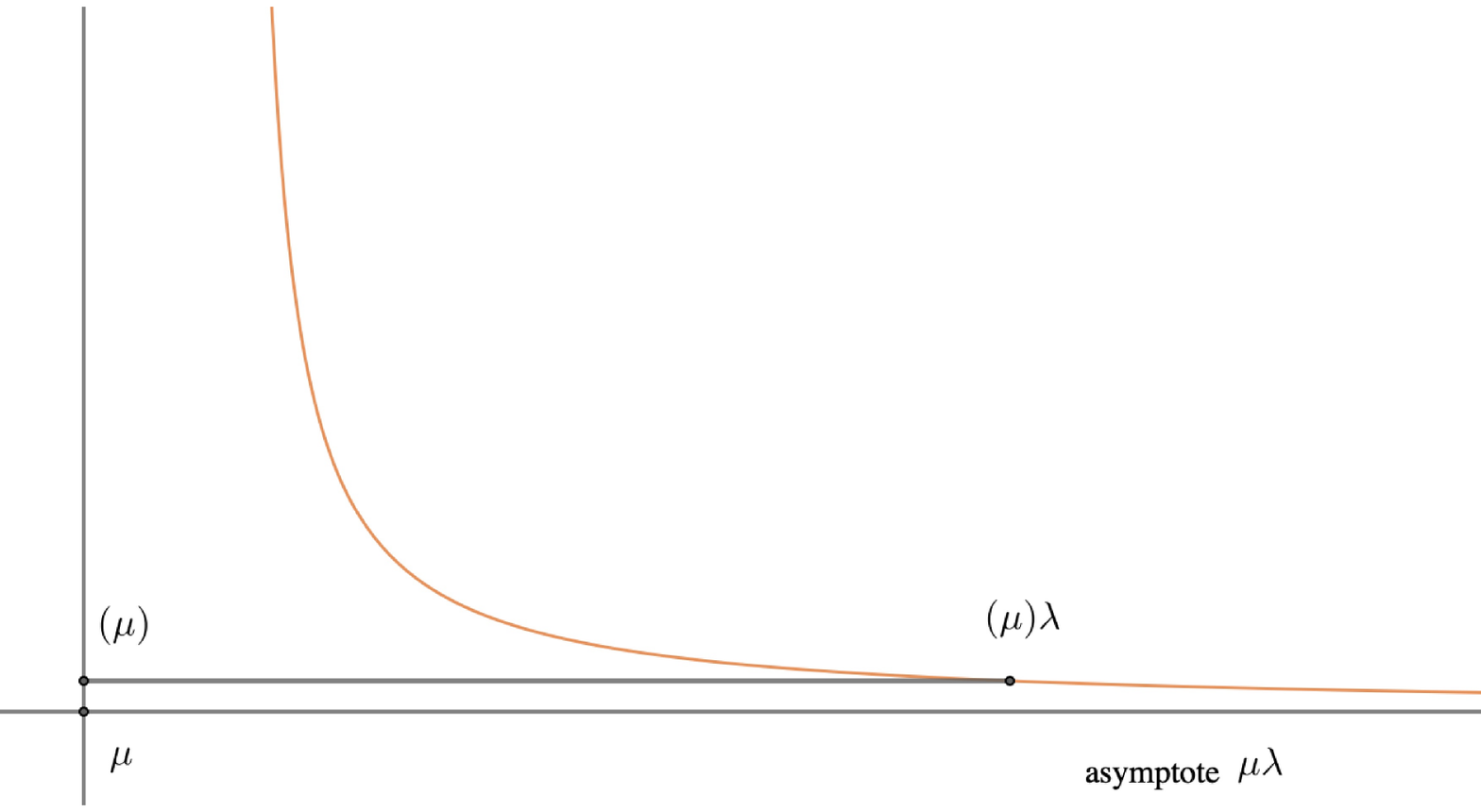}    
\end{center}
\caption{Infinitesimal $(\mu)$ and the corresponding bounded infinity
  $(\mu)\lambda$.}
\label{fig1}
\end{figure}

\subsection{Scholium}
\label{s12}

In the Scholium to Proposition\;11, Leibniz points out a further
difference between a bounded infinity and an unbounded infinity:
\begin{enumerate}\item[]
%
%
Therefore one cannot say that a bounded line is a geometric mean
between a point, which is the Minimum line, and an unbounded line,
which is the Maximum line.  But one can say that a finite line is the
geometric mean, in a sense not approximate but precise, between an
infinitely small line and an infinite line [i.e., bounded infinity].%
\footnote{``Hinc dici non potest Lineam teminatam esse proportione
  mediam inter punctum seu lineam minimam, et interminatam seu lineam
  maximam.  At dici potest lineam finitam esse mediam propotione, non
  quodammodo, sed vere exactque inter quandam infinite parvam et
  quandam infinitam'' \cite[p.\;549]{Le2012}.}
\end{enumerate}
Thus, the geometric mean of an infinitesimal and a bounded infinity
can turn out to be a finite quantity.  Infinitesimals and bounded
infinities satisfy the usual rules of arithmetic, which is not the
case for unbounded infinity.

Leibniz holds that (unlike unbounded infinity) bounded infinity is
useful in geometry and calculus, and repeatedly describes it as a
\emph{fiction}.  Both the nature of bounded infinity and the exact
meaning of its fictionality are subject to current debate among
Leibniz scholars; see e.g., Rabouin and Arthur (\cite{Ra20}, 2020),
Samuel Henry Eklund (\cite{Ek20}, 2020), Esquisabel and Raffo Quintana
(\cite{Es21}, 2021).

\section{Intersection point interpreted}
\label{s2}

Several scholars have commented on Leibnizian points at infinity,
including Knobloch, Rabouin, and Arthur.  We will examine Knobloch's
position in Section~\ref{s21}, and that of Rabouin and Arthur, in
Section~\ref{s23}.

\subsection{Hyperbola and its asymptote}
\label{s21}

Knobloch comments on Leibniz's analysis in DQA of a hyperbola~$c$ and
its asymptote~$d$ as follows:%
\footnote{Knobloch \cite[p.\;48]{Kn90}.}

\bigskip
{\small Chaque hyperbole d'un degr\'e n'importe quel [sic] a tous les
  deux axes comme asymptotes.  Ils ne concourent point ou bien ils ne
  concourent qu'apr\`es un intervalle infiniment long avec la courbe.
\[
\text{Donc on doit constater:~$d$ asymptote }
\begin{cases} 
\text{ou bien } d\cap c=\emptyset \text{ (th\'eor\`eme 11,} 
\\ 
\text{scholie: jamais; th\'eor\`eme 23: en}
\\
\text{aucun lieu)}
\\
\\
\text{ou bien } d\cap c\not=\emptyset \text{ (th\'eor\`eme 45)}
\end{cases}
\]
}

\medskip\noindent Knobloch goes on to describe the first case ($d\cap
c=\emptyset$) as occurring for the \emph{infinitum interminatum}, and
the second case ($d\cap c\not=\emptyset$), for \emph{infinitum
  terminatum}:
\begin{enumerate}\item[]
{\small N\'eanmoins ce n'est pas une contradiction parce que d'apr\`es
  l'expli\-ca\-tion leibnizienne l'infini peut \^etre
  \emph{intermin\'e ou termin\'e.}  Si l'on a en vue la premi\`ere
  possibilit\'e, on obtient l'asymptote parfaite qui ne concourt pas
  avec la courbe {\ldots} Si l'on a en vue la deuxi\`eme
  possibilit\'e, on obtient un point commun entre la droite et la
  courbe apr\`es un intervalle infiniment long.  Leibniz emploie les
  deux conceptions.%
\footnote{Knobloch \cite[p.\;48]{Kn90}; emphasis added.}
}
\end{enumerate}
Thus, Knobloch claims that the hyperbola and its asymptote have
nonempty intersection when the asymptote is a bounded infinity.  He
bases his claim upon DQA, theorem 45 (mentioned in his second case).
The relevant passage from theorem 45 is analyzed in Section~\ref{s3}
below.  This is the unique piece of evidence presented by Knobloch in
favor of his hypothesis of identification of bounded infinity and
perspective point at infinity in Leibniz.
 
Such a claim of nonempty intersection is repeated in 1994.%
\footnote{Knobloch \cite[p.\;266]{Kn94}.}
Here Knobloch writes, reasonably enough, that
\begin{enumerate}\item[]
In order to avoid contradictions we have to try to understand the
explanations of an author in their contexts. I would like to try to
demonstrate that Leibniz dealt with the infinite in a consistent
manner, although the contrary seems to be the case.  Let us consider
three examples, etc.%
\footnote{Op. cit., p.\;265.}
\end{enumerate}
Of Knobloch's three examples of an apparent contradiction, the first
is identical to his pair ``$d\cap c = \emptyset$,~$d\cap c \not=
\emptyset$'' analyzed four years earlier.%
\footnote{In Knobloch \cite[p.\;48]{Kn90}.}
Knobloch concludes that bounded infinity is a fictional entity:
\begin{enumerate}\item[]
A bounded infinite quantity is a fictitious quantity on which we rely
if we measure infinitely long but finite spaces.%
\footnote{Knobloch \cite[p.\;267]{Kn94}.}
\end{enumerate}
One finds a similar comment in 1999:
\begin{enumerate}\item[]
[Leibniz] assumed a fictive boundary point on a straight halfline
which is infinitely distant from the beginning: a bounded infinite
straight line is a fictitious quantity.%
\footnote{Knobloch \cite[p.\;97]{Kn99}.}
\end{enumerate}

\subsection{Reliability and recourse to fictionalism}
\label{s22}

However, it is unclear how describing the \emph{infinita terminata} as
fictions could explain their reliability in mathematical reasoning.
Jesseph voices a similar concern in the following terms:
\begin{enumerate}\item[]
[T]he recourse to fictionalism is insufficient on its own to make a
demonstration employing such fictions truly rigorous or convincing.%
\footnote{Jesseph \cite[p.\,196]{Je15}.}
\end{enumerate}
From the viewpoint of Knobloch's interpretation, it is difficult to
develop a coherent reading of the primary sources in Leibniz with
regard to bounded infinities and perspective points at infinity.

\subsection{Infinite quantities and ideal points}
\label{s23}

Similarly to Knobloch, Rabouin and Arthur blend infinite quantities of
the Leibnizian infinitesimal calculus and the perspective points at
infinity \`a la Desargues and Pascal.  They write:
\begin{enumerate}\item[]
[E]ven a ``clear and distinct'' concept---i.e.\;one for which we can
provide a nominal definition (allowing us to distinguish the entity in
question) may harbour a hidden contradiction, which appears when
analysing all of its constituents.  Still, one can use such concepts
for deriving truth.  This is the case with the notion of a
mathematical fiction applied to an \emph{infinite quantity}.%
\footnote{Rabouin and Arthur \cite[p.\;406]{Ra20}; emphasis added.}
\end{enumerate}
They elaborate on infinity in relation to the point at infinity of
Desargues and Pascal as follows:
\begin{enumerate}\item[]
The parallel with a point at infinity may be recalled here since this
is a notion which produces a contradiction when inserted in some
proofs of Euclid's Elements (such as I, 27, where we assume that
parallel lines meet), but which is also useful (when accompanied with
suitable demonstrations) in order to produce general geometrical
truths, such as the ones promoted by Desargues and Pascal.%
\footnote{Op. cit., p.\;407, note 14.}
\end{enumerate}
Rabouin and Arthur go on to provide the following analysis of
Leibnizian infinitesimals:
\begin{enumerate}\item[]
Thus, when Leibniz says that he understands the infinitely small to be
a fiction, this is not a way of deflecting criticism by simply
abjuring infinitesimals, as is sometimes assumed.  It means that even
though its concept may contain a contradiction, it can nevertheless be
used to discover truths, provided a demonstration can (in principle)
be given to show that its being used according to some definite rules
will avoid contradiction.%
\footnote{Op. cit., p.\;407.}
\end{enumerate}
They reiterate the claimed connection between the infinitely small and
the perspective point at infinity:
\begin{enumerate}
\item
``All of this is crucial for a proper reading of Prop.~8, the one in
  which the fiction of `infinitely small' entities will be used for
  the first time.''%
\footnote{Op. cit., p.\;418.  In fact, the \emph{infinite parvae} are
  briefly mentioned in Proposition~8
%
%
but are not used.  Instead, Leibniz gives an exhaustion argument, and
concludes in the \emph{Scholium}: 
%
%
``I went into all these details to allow geometers who encounter a
similar reasoning to avoid engaging with it, without however running
the least risk'' \cite[p.\;549]{Le2012} (translation ours).  Leibniz
does not introduce the notation~$(\mu)\mu$ for an infinitely small
line until Proposition~11.  }
\item
``The idea of fiction is mentioned a first time for designating the
  `point at infinity' introduced by geometers developing projective
  considerations, such as Desargues and Pascal (schol.\;VII).''%
\footnote{Rabouin and Arthur \cite[p.\;418, note 42]{Ra20}.}
\end{enumerate}
The alleged connection between infinitesimals and Desargues is
mentioned yet again:
\begin{enumerate}\item[]
The parallel between the introduction of point at infinite distance
and infinitesimals {\ldots}\;appears in the DQA and was already
present in Desargues, {\ldots} \cite[p.\;421, note 52]{Ra20}
\end{enumerate}

A debate of long standing concerns the issue of the fictionality of
the Leibnizian infinitesimal and infinite quantities.  In his articles
\cite{Kn90} and~\cite{Kn94}, Knobloch seeks to relate such
fictionality to alleged paradoxical behavior of bounded infinity when
the hyperbola and its asymptote are said to meet, while Rabouin and
Arthur read the fictionality of infinitesimals as inconsistency in
\cite{Ra20}.  In Sections~\ref{s3} and \ref{s4} we will analyze the
interpretation of fictionality and of the ideal intersection point
between a curve and its asymptote.

\section{Theorem 45 of \emph{De Quadratura Arithmetica}}
\label{s3}

\begin{figure}
\begin{center}
\includegraphics
[scale=0.5,
trim = 0 0 0 0,clip]
{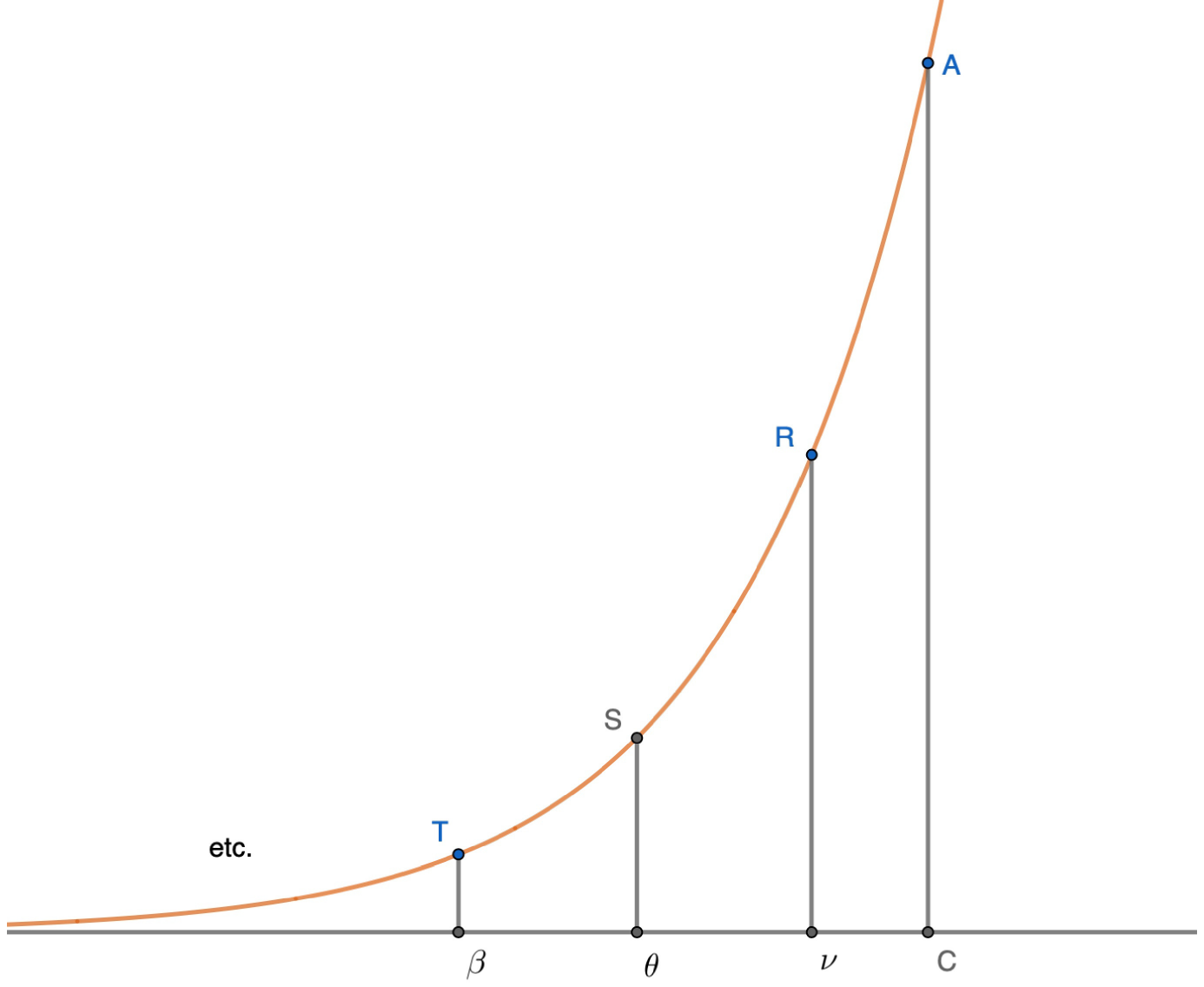}    
\end{center}
\caption{Curve $ARST$ and asymptote $C\beta$.} 
\label{fig2}
\end{figure}

Knobloch's interpretation relies on a reading of a sentence in the
proof of theorem 45 (i.e., \emph{Propositio} XLV) of DQA.\, The
sentence reads as follows:
\begin{enumerate}\item[]
%
%
It remains to show that the line $C\beta$\;etc.  represents an
asymptote, i.e., that it cannot meet the logarithmic curve
$ARST$\;etc.  except at an infinite distance.%
\footnote{``Superest tantum ut ostendamus rectam $C\beta$\;etc.  esse
  asymptoton seu non nisi in infinito abhinc intervallo occurrere
  posse curvae Logaritmicae, $ARST$\;etc.''  \cite[p.\;633]{Le2012}.}
\end{enumerate}
The issue is the nature of the intersection between the asymptote
denoted ``$C\beta$\;etc.''  and ``$ARST$\;etc.'' denoting the
logarithmic curve.  Leibniz supplements the notation for a line by the
abbreviation ``etc."  to indicate that he is referring to an unbounded
line rather than a bounded infinity.  Thus, he uses the notation
``$C\beta$\;etc.''  three times in the proof of theorem 45 to denote
the unbounded asymptote.
%
%
The figure illustrating the theorem%
\footnote{See Figure 14 in \cite[p.\;624]{Le2012}.}
%
%
marks the finite points~$C$ and~$\beta$ lying on the asymptote to the
curve; see Figure~\ref{fig2}.

We note the following points:
\begin{itemize}
\item 
The term \emph{bounded infinity} is not used in the proof of
theorem~45.
\item 
The intersection at infinity between the logarithmic curve
denoted~$ARST$\;etc.  and its asymptote
is only mentioned in passing in a double negation:%
\footnote{A similar double negation referring to the intersection at
  infinity of a \emph{hyperbola} and its asymptote occurs in the proof
  of Proposition 22.
}
``cannot meet {\ldots}~except at an infinite distance''.%
\footnote{Knobloch's paraphrase of the phrase transforms the double
  negative clause into a positive one: ``Second assertion: The
  straight line~$C\beta$\;etc.  is an asymptote, which is
  (\emph{seu}), it can meet (\emph{occurrere posse}) the logarithmic
  curve~$ARST$ only after an infinite interval (\emph{infinito abhinc
    intervallo})'' \cite[p.\;28]{Kn18}.}
\end{itemize}
We can therefore make the following five remarks.

\medskip\noindent 
1. If one wishes to explain why, according to Leibniz, bounded
infinity is useful in geometry while unbounded infinity is not,
theorem 45, for all its intrinsic interest, is of little help.

\medskip\noindent 
2. In Proposition 11, \emph{both} the unbounded asymptote and a
genuine bounded infinite~$(\mu)\lambda$ occur in the same proof, and
Leibniz repeatedly states that they are distinct.

\medskip\noindent 
3. Using the terminology of Proposition 11, Knobloch acknowledges the
existence of an infinitely thin rectangle formed by drawing a pair of
lines, parallel to the axes, and passing through a point infinitely
far along the curve~$\lambda$ with infinitesimal abscissa: ``Si
l'abscisse~$\mu(\mu)$ est infiniment petite,
l'ordonn\'e~$(\mu)\lambda$ est infiniment longue, c'est \`a dire plus
grande qu'une droite quelconque designable, le rectangle effectu\'e
par la droite infinie et la droite infiniment petite est \'egal \`a un
carr\'e fini constant selon la nature de l'hyperbole.''%
\footnote{Knobloch \cite[p.\;41]{Kn90}.}
%
%
The side~$(\mu)\lambda$ is an instance of a bounded infinity.  It is a
\emph{subline} of the unbounded infinite line which passes through the
infinitesimal abscissa~$(\mu)$ and is parallel to the asymptote.  By
the nature of such a rectangle (infinitely thin, infinitely long),
there is always a \emph{nonzero} distance from the vertex (i.e.,
rectangle vertex lying on the curve) to the asymptote, even when the
point is infinitely close to the asymptote because~$(\mu)$ is
infinitesimal.

\medskip\noindent 
4.  Leibniz observes that bounded lines (i.e., segments) \emph{cannot
exhaust} the unbounded infinite line.%
\footnote{The full observation reads as follows: 
%
%
``Indeed, just as one does not change a bounded line by adding to or
  removing from it some points, even infinitely many, so also by
  reproducing any number of times a bounded line, one can neither
  constitute nor exhaust an unbounded line.  It is otherwise for a
  bounded infinite line, which can be conceived as formed of a
  multitude of finite lines, even though such a multitude exceeds all
  number'' \cite[p.\;549]{Le2012} (translation ours).}
The claim that they cannot exhaust it presupposes in particular that
it is meaningful to envision an attempt to carry out such an
exhaustion; namely that a bounded line, while unable to exhaust it, is
a \emph{subline} of the unbounded infinite line (as in the example
$(\mu)\lambda$ mentioned in item 3 above).  If so, how could the
unbounded infinite line~$C\beta$\;etc.  have empty intersection with
the logarithmic curve~$ARST$\;etc.  (as claimed by Knobloch; see
Section~\ref{s2}) whereas its subline nonetheless manages to meet this
curve?

\medskip\noindent 
5.  Leibniz's \emph{infinita terminata} can be naturally scaled; they
can be multiplied by~$2, 3, \ldots$ just as the infinitesimal~$(\mu)$
can be divided by~$2,3,\ldots$, and satisfy the usual rules of
arithmetic (see Section~\ref{s12}).  One cannot obtain such properties
by adjoining a \emph{single} point at infinity, \`a la projective
geometry, implied in Leibniz's comment quoted in Section~\ref{s2}.  It
is unclear how calling it ``fictional'' could help here (see Jesseph's
comment in Section~\ref{s22}).

\medskip 
It is therefore difficult to develop a consistent reading of the
Leibnizian notions of bounded infinity and perspective point at
infinity from the viewpoint of Knobloch's interpretation that seeks to
identify them.

Projective geometry was only in its incipient stages at the time, but
Leibniz did speak of ideal meeting points at infinity for parallel
lines (see Section~\ref{s4}).  We argue that, if the hyperbola and its
asymptote meet at infinity, it is only in the sense familiar from
projective geometry, a sense distinct from the \emph{infinita
  terminata} of Leibniz's geometry and calculus.  We review the
history of the idea of the perspective point at infinity in
Section~\ref{s3b}.

\section{Infinite distance from Kepler to Desargues}
\label{s3b}

By the last quarter of the 17th century when Leibniz started work on
DQA, the idea of an infinitely distant point in the geometry of
perspective was already a familiar one, due to the work of Kepler,
Desargues, and Bosses.

In 1604, Kepler referred to an infinitely distant point associated
with a pencil of parallel lines as a ``blind focus.''  He held that
\begin{enumerate}
\item[] In the Parabola one focus,~$D$, is inside the conic section,
  the other is to be imagined either inside or outside, lying on the
  axis [of the curve] at an infinite distance from the former
  (\emph{alter vel extra vel intra sectionem in axe fingendus est
    infinito intervallo \`a priore remotus}), so that if we draw the
  straight line~$HG$ or~$IG$ from this blind focus (\emph{ex illo
    caeco foco}) to any point~$G$ on the conic section, the line will
  be parallel to the axis~$DK''$ (\emph{Ad Vitellionem Paralipomena},
  cap.\;IV n.\;4: \emph{De coni sectionibus}).%
\footnote{For more details see \cite[p.\;450]{Fi86},
  \cite[pp.\,186--187]{Fi87}, and \cite[p.\;568]{De16}.}
\end{enumerate}

According to Debuiche, the idea that the projective closure of a line
is a circle can already be detected in Desargues \cite{De39}:
\begin{enumerate}
\item[] Desargues presents the idea of a complete correspondence
  between a straight line and a circle, since a straight line can be
  considered as a circle closed in on itself at an infinite distance.%
\footnote{Debuiche \cite[p.\;373]{De13}.}
\end{enumerate}
Field and Gray note: 
\begin{enumerate}\item[]
Desargues began with the remark that lines will be supposed to contain
a point at infinity, which may be reached by travelling in either
direction along the line.%
\footnote{Field and Gray  \cite[p.\;47]{Fi87}.}
\end{enumerate}
In connection with Leibniz's work on perspective, Rodr\'{\i}guez
Hurtado et al.~note that
\begin{enumerate}
\item[] The characterisation of the point of view as the meeting point
  in the infinity of the parallels comes from Arguesian
  perspective.$^{56}$ \cite[p.\,16]{Ro21}%
\footnote{Their note 56 reads: ``$^{56}$It is worth mentioning that
  Leibniz transcribes the definition of point of view that Desargues
  makes at the end of the \emph{Brouillon Project} (included in: A VII
  7, 111).''  The page number given is incorrect.  It should be
  A~VII~7, item 65, page 593.}
\end{enumerate}
The reference is to Desargues' 1639 \emph{Brouillon Project}
\cite{De39}.  Based on Leibniz's 1679 letter to Huygens,
Rodr\'{\i}guez Hurtado et al.~point out a significant difference
between the Cartesian model of magnitudes and Leibniz's search for a
qualitative model for studying perspective:
\begin{enumerate}\item[]
In contrast to the Cartesian algebraic model, centred on the
determination of magnitudes (analysis of quantitative relations),
Leibniz wanted to construct a qualitative model, based on analysis of
position (\emph{situm}).%
\footnote{Rodr\'{\i}guez et al.~\cite[p.\;4]{Ro21}.}
\end{enumerate}
While critical of certain aspects of \cite{Ro21}, a recent text by
Debuiche and Brancato confirms that ``perspective science can be
understood as containing the `whole \emph{Geometria Situs}' since the
geometry of situation only deals with the mutual positions of points
in space.''%
\footnote{Debuiche and Brancato \cite[p.\;67]{De23}.}
Accordingly, perspective points at infinity are not expected to have
the properties of magnitudes, unlike Leibnizian \emph{infinita
  terminata}.

Leibniz explicitly acknowledged a debt to Desargues in a 1692
publication in \emph{Act.\;Erudit.\;Lips.} entitled \emph{De Linea ex
  Lineis Numero Infinitis}:

\begin{enumerate}\item[] 
%
%
Geometers customarily refer as ordinates to parallel lines in any
number, traced between a curve and a fixed line (directrix); when they
are perpendicular to the latter (which then plays the role of an
axis), one refers to ordinates par excellence.  Desargues generalized
this by considering also as ordinates, the lines which converge toward
a unique common point or diverge from it.  Parallel lines can be
considered converging or diverging lines, if one fictively considers
that their common point is at an infinite distance.%
\footnote{``\emph{Ordinatim applicatas} vocare solent Geometrae rectas
  quotcunque inter se parallelas, quae a curva ad rectam quandam
  (directricem) usque ducuntur, quae cum ad directricem (tamquam axem)
  sunt normales, solent vocari \emph{ordinatae} $\kappa\alpha\tau$'
  '\!\!$\varepsilon\xi o\chi\acute{\eta}\nu$.
%
%
Desarguesius rem prolatavit, et sub ordinatim applicatis etiam
comprehendit rectas \emph{convergentes} ad unum punctum commune, aut
ab eo \emph{divergentes}.  Et sane parallelae sub convergentibus aut
divergentibus comprehendi possunt, fingendo punctum concursus infinite
ab hinc distare'' (see Gerhardt \cite{Ge50}, vol.\;V, pp.\;266--267).}
\end{enumerate}
Leibniz deepened his knowledge of projective methods in Hannover by
studying Abraham Bosse's writings on Desargues as noted in
\cite{Ro21}, but the more detailed statement in 1692 is merely a
clearer elaboration of the position already found in DQA.

In Section~\ref{s4}, we will analyze the occurrence in Leibniz of the
Arguesian point at infinity to which parallel lines `diverge.'

\section{Leibniz's pencil of parallel lines}
\label{s4}

In the \emph{Scholium} to Proposition 7 of DQA, Leibniz speaks of a
pencil of parallel lines meeting at a fictional point~$A$.  We argue
that this is a fiction in a sense distinct from his infinitesimals and
\emph{infinita terminata}.  Leibniz attributes the idea to
``illustrious geometers'': 
\begin{enumerate}\item[]
%
Furthermore, illustrious geometers having undertaken to study the
conics from a general viewpoint, call ordinates to curves not only, as
is common, the parallel lines ${}_1C\,_1B$, ${}_2C\,_2B$,
${}_3C\,_3B$, but also the lines $A\,_1C$, $A\,_2C$,~$A\,_3C$ which
all converge toward a unique point $A$ (which is entirely correct
since one can, without committing an error, consider parallel lines as
convergent lines, up to considering fictively [\emph{fingatur}] that
their point of intersection or their common center is at an infinite
distance, as in the case of the focus or vertex of the parabola.%
\footnote{``Porro cum Clarissimi Geometrae, qui Conica universaliter
  tractare coepere, ordinatarum ad curvas nomine comprehendant non
  tantum rectas parallelas, quales sunt ${}_1C\,{}1B$, ${}_2C\,_2B$,
  ${}_3C\,_3B$, ut vulgo fieri solet, sed etiam rectas $A\,_1C$,
  $A\,_2C$, $A\,_3C$, quae omnes ad unum punctum commune $A$,
  convergunt (quod vel ideo recte fit, quoniam ipsaemet 5 parallelae
  sine errore pro convergentibus sumi possunt, ita tantum ut punctum
  concursus earum seu centrum commune infinite abesse fingatur,
  quemadmodum alter parabolae focus aut vertex)''
  \cite[p.\;538]{Le2012}.  Leibniz's reference to a parabola can be
  interpreted as follows.  If one sends rays out of the (finite) focal
  point of the parabola, then after bouncing off the parabola, the
  rays turn into a pencil of parallel lines.  The corresponding ideal
  point at infinity~$A$ is the ``focal point at infinity'' of the
  parabola.  See Section~\ref{s3b} for a discussion of possible
  antecedents in Kepler and Desargues.}

\end{enumerate}
Parmentier identifies the illustrious geometers as Pascal and
Desargues.%
\footnote{Parmentier \cite[p.\;73]{Le04b}.}

A curve (such as a hyperbola or a logarithmic curve) and its asymptote
are obviously not a pair of parallel lines, but they meet at infinity
at a unique fictional point~$A$ determined by the pencil of lines
parallel to the asymptote.

To illustrate that the perspective point at infinity can be easily
formalized as the ideal point of intersection between the logarithmic
curve and its asymptote, we outline a modern formalisation in the case
of the hyperbola \mbox{$xy=1$} and its horizontal asymptote~$y=0$ in
the affine plane.  Passing to homogeneous coordinates~$[x_1,x_2,x_3]$
where~$x=\frac{x_1}{x_3}$ and~$y=\frac{x_2}{x_3}$, we obtain the
equation
\begin{equation}
\label{e51}
x_1 x_2 = x_3^2
\end{equation}
for the projective completion of the hyperbola, and equation~$x_2=0$
for that of the asymptote.  The point at infinity for the asymptote is
the point~$A=[1,0,0]$.  The point~$A$ clearly lies on the
curve~$x_1x_2=x_3^2$ of equation~\eqref{e51}, as well, and can
therefore be thought of as the ideal point of intersection at infinity
between the hyperbola and its horizontal asymptote.

What could be the relation between such an ideal point at infinity~$A$
and Leibniz's bounded infinities?  We make the following two
observations.

\begin{enumerate}
\item
Bounded infinity~$(\mu)\lambda$, as well as the infinitesimal
$(\mu)\mu$, are only discussed by Leibniz in Proposition\;11.
Therefore the mention of the ideal point~$A$ earlier in the text,
namely in the \emph{Scholium} following Proposition 7, could not be
related to bounded infinity, unless we presume Leibniz to be sloppy in
his exposition in DQA by using a concept before discussing it.

\item
The ideal point at infinity~$A$ in projective geometry does not make
sense as an infinite \emph{magnitude}, because it cannot occur as an
element in an ordered system that is greater than all the other
magnitudes.  Indeed, such a point~$A$ is assigned to a pencil of
\emph{unoriented} (undirected) lines; if~$A$ were greater than all the
other magnitudes, it would also have to be declared smaller than all
the other magnitudes, leading to an absurdity.%
\footnote{Thus if one starts with an affine line~$\R$, the
  corresponding projective line~$\R\mathbb P^1$ is the circle~$S^1$
  (which admits no natural structure of an order).  See
  Section~\ref{s3b} for a discussion of possible antecedents in
  Desargues.  Anglade and Briend note that ``Desargues \'enonce
  {\ldots}~une analogie entre le cercle et la droite, ce qui peut
  laisser croire qu'il avait une image juste de ce que l'on appelle
  aujourd'hui la droite projective (r\'eelle) comme \'etant
  topologiquement un cercle'' \cite[p.\;550]{An17}.  Anglade and
  Briend published a series of in-depth studies of the \emph{Brouillon
    Project} culminating in \cite{An22}.}
\end{enumerate}
Thus, Leibniz's (projective) ideal point at infinity~$A$ cannot be
identified with an \emph{infinitum terminatum} without creating
unnecessary inconsistencies.  An \emph{infinitum terminatum}, being
the inverse of an infinitesimal, cannot be a perspective ideal point
at infinity.

There are hints in Knobloch's later work that he is aware of the
difficulties with his interpretation of Leibnizian infinitesimals, as
when he writes:
\begin{enumerate}\item[] 
[T]he error is smaller than any assignable error and therefore zero
{\ldots}~Such an error necessarily is equal to zero as Leibniz rightly
states.  For if we assume that such an error is unequal to zero it
would have a certain value.  But this implies a contradiction against
the postulate that the error has to be smaller than any assignable
quantity, that is, also smaller than this certain value.  Yet, Leibniz
explicitly calls such errors infinitely small: \emph{We should not try
  to make things seem better}.  \cite[p.\;12]{Kn18} (emphasis added)
\end{enumerate}
Leibnizian infinitesimals were not assignable, as made clear by
Leibniz himself, who wrote: 
\begin{enumerate}\item[]
Bien qu'elles ne soient pas assignables, elles se trouvent \^etre
quelque chose d’existant et non pas un rien absolu.%
\footnote{Leibniz as quoted by Bella in \cite[p.\,195]{Be19}.}
\end{enumerate}

As analyzed above, there is no contradiction between the two theorems
(Propositio 11 and Propositio 45), as they deal with different notions
of infinity.  Propositio 11 deals with the \emph{infinitum
  terminatum}, whereas Propositio 45 deals with a perspective point at
infinity.  In a recent text, Knobloch discusses both theorems and
concludes as follows: ``Here [i.e., in Theorem 45], Leibniz says the
opposite of what he had said in the demonstration of theorem 11.
There, he had explicitly excluded the possibility that the curve meets
the asymptote in accordance with the meaning of its name''
\cite[p.\;9]{Kn24}.  Thus, Knobloch prefers to attribute a
contradiction to Leibniz rather than accept an interpretation other
than the syncategorematic one.  Knobloch adds the following sentence:
``We shall come back to this matter of fact'' (ibid.), but he never
does.


\section{Conclusion}

We have sought to redress a conflation of Leibniz's notion of bounded
infinity and the notion of the perspective (projective) point at
infinity, in the recent literature on Leibniz.  We argued that
Leibniz's proof of theorem\;45 of his \emph{De Quadratura Arithmetica}
provides no evidence for relating his notion of bounded infinity to a
remote point of intersection (ideal projective point) of a curve
(hyperbola or logarithmic curve) and its asymptote.  Furthermore, the
hypothesis of such a point of intersection at \emph{bounded infinite}
range is mathematically incoherent.  

Postulating that there are only two types of infinity in Leibniz --
bounded and unbounded -- leads to a paradoxical conclusion that what
is identifiably the projective ideal point at infinity of a pencil of
parallel lines, must be bounded infinity.  Such an approach leads to
unnecessary inconsistencies (as detailed in Section~\ref{s3} and
\ref{s4}).  Conflating perspective points at infinity and
\emph{infinita terminata} amounts to what Leibniz may have called a
``category error'': the former belong in analysis situs whereas the
latter belong in analysis of magnitudes.

As noted by Jesseph (see Section~\ref{s22}), appeals to fictionality
are insufficient to provide an adequate account of Leibniz's
mathematical procedures.  More convincing accounts of the fictionality
of infinitesimals were recently developed e.g., in \cite{Ek20},
\cite{Es21}.  Interpreting Leibniz’s infinitesimals is an area of
lively debate.  In 2021, Bair et al.~published a comparative study of
three interpretations \cite{21a}.  Katz et al.~\cite{21g} presented
three case studies in Leibniz scholarship.  In 2022, Katz et
al.~presented and analyzed a pair of rival approaches \cite{22b}.  In
the same year, Archibald et al.~formulated some criticisms
\cite{Ar22b}.  In 2023, Bair et al.~provided both a brief response
\cite{23a} and a detailed response \cite{22a}.

In sum, one is led to recognize that there are multiple types of
infinity in Leibniz's geometry and calculus.  We conclude that an
ideal point at infinity (associated with a pencil of parallel lines)
is borrowed from Desargues, and is a type distinct from the
\emph{infinita terminata}.  Such an approach enables a consistent
interpretation of the Leibnizian notions of bounded infinity and
perspective point at infinity, and explains how the former (but not
the latter) can satisfy the usual rules of arithmetic such as
invertibility and scalability.

\end{document}